\documentclass{article}

\usepackage{arxiv}

\usepackage[utf8]{inputenc} % allow utf-8 input
\usepackage[T1]{fontenc}    % use 8-bit T1 fonts
\usepackage{hyperref}       % hyperlinks
\usepackage{url}            % simple URL typesetting
\usepackage{booktabs}       % professional-quality tables
\usepackage{amsfonts}       % blackboard math symbols
\usepackage{nicefrac}       % compact symbols for 1/2, etc.
\usepackage{microtype}      % microtypography
\usepackage{lipsum}         % Can be removed after putting your text content
\usepackage{graphicx}
\usepackage{doi}

\usepackage{amssymb,amsthm,amsmath}
\usepackage{xcolor,paralist,hyperref,titlesec,fancyhdr,etoolbox}
\newtheorem{theorem}{Theorem}[]
\newtheorem{definition}[theorem]{Definition}

\newtheorem*{remark}{Remark}

\title{On the existence of directed strongly regular graphs with parameters (22, 9, 6, 3, 4)\thanks{Translated from the Russian original published in: 	V. A. Byzov, I. A. Pushkarev, “On the existence of directed strongly regular graphs with parameters (22,9,6,3,4)”, Prikladnaya Diskretnaya Matematika, 2024, no. 66, 86--96}}

\date{September, 2025}
\author{ \href{https://orcid.org/0000-0002-3613-5949}{\includegraphics[scale=0.06]{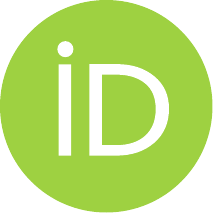}\hspace{1mm}Viktor A.~Byzov} \\
	Vyatka State University\\
	Kirov,  Russian Federation \\
	\texttt{vbyzov@yandex.ru} \\
	\And
	{ \href{https://orcid.org/0000-0002-3610-2872}{\includegraphics[scale=0.06]{orcid.pdf}\hspace{1mm}Igor A.~Pushkarev}} \\
	Vyatka State University\\
	Kirov,  Russian Federation \\
	\texttt{god\_sha@mail.ru}\\
}

\renewcommand{\headeright}{~}
\renewcommand{\undertitle}{~}
\renewcommand{\shorttitle}{On the existence of directed strongly regular graphs with parameters (22, 9, 6, 3, 4)}

\begin{document}
\maketitle

\begin{abstract}
The paper shows the existence of a family of directed strongly regular graphs with parameters \mbox{(22, 9, 6, 3, 4)}. The adjacency matrices of the found digraphs are composed of $3\times 3$ circulant blocks. The automorphism group of all the digraphs found is the group~$\mathbb{Z}_3$. The structure of the resulting digraphs is described using concepts of skeleton and rigging.
\end{abstract}

% keywords can be removed

\vspace{10pt}
\section{Introduction}

We introduce the basic concepts and notations necessary for the subsequent exposition.

In this work, we consider only directed graphs (hereafter referred to as digraphs) without loops and without multiple edges in the same direction. The vertices of such digraphs will be numbered by natural numbers starting from one. Adjacency matrices will be formed in one of the standard ways for digraphs, namely: we place a 1 in the $i$-th row and the $j$-th column if there exists an edge in the digraph going from vertex $i$ to vertex $j$ (with rows and columns also numbered by natural numbers). All other entries of the matrix are set to 0.

For the identity matrix of order $n$, we use the notation $I_n$; for the square matrix of order $n$ consisting entirely of ones, we use the notation $J_n$. The notation $J_{k,l}$ denotes a $k \times l$ matrix consisting of ones. The segment of natural numbers ${1, 2, \ldots, n}$ will be denoted by $[n]$.

The concept of a strongly regular digraph was introduced by A.\,M.~Duval in~\cite{bibDuval1} as an oriented generalization of the concept of strongly regular undirected graphs (see~\cite{bibSRG1}). We now present two equivalent definitions of such a digraph.

\begin{definition}
    \label{def:main1}
    A strongly regular digraph with parameter set $(v, k, t, \lambda, \mu)$ is a digraph on $v$ vertices satisfying the following conditions:
    \begin{enumerate}
        \item the outdegree and indegree of every vertex are equal to $k$;
        \item for any vertex $x$, there are exactly $t$ paths of the form $x \to z \to x$;
        \item for any edge $x \to y$, there are exactly $\lambda$ paths of the form $x \to z \to y$;
        \item if there is no edge $x \to y$ in the digraph, then there are exactly $\mu$ paths of the form $x \to z \to y$.
    \end{enumerate}
\end{definition}

\begin{definition}
    \label{def:main2}
    A strongly regular digraph with parameter set $(v, k, t, \lambda, \mu)$ is defined as a digraph whose adjacency matrix $A$ satisfies the following relations:
    \begin{equation}
        \label{eq:main1}
        A^2 = t I_v + \lambda A + \mu (J_v - I_v - A),
    \end{equation}
    \begin{equation}
        \label{eq:main2}
        A J_v = J_v A = k J_v.
    \end{equation}
\end{definition}

Digraphs of the described type will be denoted as $\text{dsrg}(v, k, t, \lambda, \mu)$. In~\cite{bibDuval1}, a number of necessary conditions are given that the parameters of a strongly regular digraph must satisfy. However, as in the case of strongly regular graphs, there are many parameter sets for which it is unknown whether a strongly regular digraph with these parameters exists. A.\,E.~Brouwer, on the website~\cite{bibBrouwer1}, systematizes information about which parameter sets admit such digraphs (and provides these digraphs) and for which sets the question of existence remains open. Up to the present moment, the smallest such case (in terms of the number of vertices) was the existence question for the digraph $\text{dsrg}(22, 9, 6, 3, 4)$. In this work, we construct a family of digraphs with this parameter set.

\begin{remark}
    We note a certain inconsistency in the order of listing the parameters of a strongly regular digraph among different authors. In the foundational paper~\cite{bibDuval1}, the following parameter order is used: $(v, k, \mu, \lambda, t)$. In other works (for example,~\cite{bibBrouwer631}), the parameters are listed in the order $(v, k, t, \lambda, \mu)$. The second method appears more natural, since it aligns better with the standard notation for strongly regular graphs $\text{srg}(v, k, \lambda, \mu)$. Therefore, in this work we adopt this latter convention.
\end{remark}

We will need the concept of a circulant matrix. A circulant matrix (or circulant) is a square matrix of the form
\begin{equation}
\label{eq:circul}
\begin{pmatrix}
a_1 & a_2 & \dots & a_n \\
a_n & a_1 & \dots & a_{n-1} \\
\vdots & \vdots & \ddots & \vdots \\
a_2 & a_3 & \dots & a_1
\end{pmatrix},
\end{equation}
that is, a matrix in which each row, starting from the second, is obtained by a cyclic shift of the previous row one position to the right.

It is true that the set of circulant matrices of order $m$ with elements from $\mathbb{Z}$ (under the standard operations of multiplication and addition) forms a ring isomorphic to the factor ring $\mathbb{Z}[x]/(x^m - 1)$ (see, for example, \cite{bibKra1}). Under this isomorphism, the matrix \eqref{eq:circul} corresponds to the polynomial $a_1 + a_2 x + \ldots + a_n x^{n-1}$.

We call the \textit{compactification} of a block matrix $M$ consisting of circulants the process of replacing all circulants by their corresponding polynomials under the isomorphism described above. The resulting matrix will be denoted by $M(x)$. It is easy to see that the compactification of matrices is compatible with the operations of matrix addition and multiplication. That is, if square block matrices $M_1$, $M_2$, $M_3$, and $M_4$ of the same order, consisting of $m \times m$ circulants, satisfy $M_3 = M_1 + M_2$ and $M_4 = M_1 \cdot M_2$, then it holds that $M_3(x) \equiv M_1(x) + M_2(x) \pmod{x^m - 1}$ and $M_4(x) \equiv M_1(x) \cdot M_2(x) \pmod{x^m - 1}$.

As an illustration, consider the adjacency matrix $S$ of the Shrikhande graph (see~\mbox{\cite{bibSRG1, bibShrikhande1}}), which is strongly regular with parameter set $(16, 6, 2, 2)$. With a suitable numbering of the vertices, it can be represented in the following form:
\begin{equation}
S = \left(
\begin{array}{cccc|cccc|cccc|cccc}
0 & 0 & 0 & 0 & 0 & 0 & 1 & 1 & 0 & 0 & 1 & 1 & 0 & 1 & 0 & 1 \\
0 & 0 & 0 & 0 & 1 & 0 & 0 & 1 & 1 & 0 & 0 & 1 & 1 & 0 & 1 & 0 \\
0 & 0 & 0 & 0 & 1 & 1 & 0 & 0 & 1 & 1 & 0 & 0 & 0 & 1 & 0 & 1 \\
0 & 0 & 0 & 0 & 0 & 1 & 1 & 0 & 0 & 1 & 1 & 0 & 1 & 0 & 1 & 0 \\ \hline
0 & 1 & 1 & 0 & 0 & 0 & 0 & 0 & 0 & 1 & 0 & 1 & 0 & 0 & 1 & 1 \\
0 & 0 & 1 & 1 & 0 & 0 & 0 & 0 & 1 & 0 & 1 & 0 & 1 & 0 & 0 & 1 \\
1 & 0 & 0 & 1 & 0 & 0 & 0 & 0 & 0 & 1 & 0 & 1 & 1 & 1 & 0 & 0 \\
1 & 1 & 0 & 0 & 0 & 0 & 0 & 0 & 1 & 0 & 1 & 0 & 0 & 1 & 1 & 0 \\ \hline
0 & 1 & 1 & 0 & 0 & 1 & 0 & 1 & 0 & 0 & 0 & 0 & 1 & 1 & 0 & 0 \\
0 & 0 & 1 & 1 & 1 & 0 & 1 & 0 & 0 & 0 & 0 & 0 & 0 & 1 & 1 & 0 \\
1 & 0 & 0 & 1 & 0 & 1 & 0 & 1 & 0 & 0 & 0 & 0 & 0 & 0 & 1 & 1 \\
1 & 1 & 0 & 0 & 1 & 0 & 1 & 0 & 0 & 0 & 0 & 0 & 1 & 0 & 0 & 1 \\ \hline
0 & 1 & 0 & 1 & 0 & 1 & 1 & 0 & 1 & 0 & 0 & 1 & 0 & 0 & 0 & 0 \\
1 & 0 & 1 & 0 & 0 & 0 & 1 & 1 & 1 & 1 & 0 & 0 & 0 & 0 & 0 & 0 \\
0 & 1 & 0 & 1 & 1 & 0 & 0 & 1 & 0 & 1 & 1 & 0 & 0 & 0 & 0 & 0 \\
1 & 0 & 1 & 0 & 1 & 1 & 0 & 0 & 0 & 0 & 1 & 1 & 0 & 0 & 0 & 0
\end{array}
\right).
\end{equation}

The matrix $S$ can be compactified into the matrix:
\begin{equation}
S(x) =
\begin{pmatrix}
    0 & x^2+x^3 & x^2+x^3  & x+x^3\\
    x+x^2 & 0 & x+x^3 & x^2+x^3\\
    x+x^2 & x+x^3 & 0 & 1+x\\
    x+x^3 & x+x^2 & 1+x^3 & 0
\end{pmatrix}.
\end{equation}

Moreover, the original matrix satisfies the equation $S^2 = 4I_{16} + 2J_{16}$; consequently, the matrix $S(x)$ satisfies the congruence
\begin{equation}
S^2(x) \equiv 4I_4 + 2(1 + x + x^2 + x^3) J_4 \pmod{x^4 - 1}.
\end{equation}

\section{Search for digraphs $\text{dsrg}(22, 9, 6, 3, 4)$}
\label{sec1}

According to Definition~\ref{def:main2}, we need to find a binary square matrix of order 22 satisfying conditions~\eqref{eq:main1} and~\eqref{eq:main2}. A complete exhaustive search over all such matrices would take an unmanageably long time, so we will restrict the search to matrices of a special form. The idea of the approach described below is adapted from the work of O.~Gritsenko~\cite{bibGritsenko1}.

We will search for the adjacency matrix $A$ of the digraph $\text{dsrg}(22, 9, 6, 3, 4)$ in the form presented in formula~\eqref{eq:block_repr}. Note that the first column and first row are chosen so that the first two conditions of Definition~\ref{def:main1} are automatically satisfied. In these row and column, zeros replace the dots in the given representation. The matrices $K_{ij}$ ($1 \leqslant i, j \leqslant 7$) in this representation are circulants of order three.
\begin{equation}
\label{eq:block_repr}
A = \left(
\begin{array}{c|ccc|ccc|ccc|c|ccc}
  0 & 1 & 1 & 1 & 1 & 1 & 1 & 1 & 1 & 1 & \ldots & 0 & 0 & 0 \\ \hline
  1 & & & & & & & & & & & & & \\
  1 & & K_{11} & & & K_{12} & & & K_{13} & & \ldots & & K_{17} & \\
  1 & & & & & & & & & & & & & \\ \hline
  1 & & & & & & & & & & & & & \\
  1 & & K_{21} & & & K_{22} & & & K_{23} & & \ldots & & K_{27} & \\
  1 & & & & & & & & & & & & & \\ \hline
  0 & & & & & & & & & & & & & \\
  0 & & K_{31} & & & K_{32} & & & K_{33} & & \ldots & & K_{37} & \\
  0 & & & & & & & & & & & & & \\ \hline
  1 & & & & & & & & & & & & & \\
  1 & & K_{41} & & & K_{42} & & & K_{43} & & \ldots & & K_{47} & \\
  1 & & & & & & & & & & & & & \\ \hline
  \vdots & & \vdots & & & \vdots & & & \vdots & & \ddots & & \vdots & \\ \hline
  0 & & & & & & & & & & & & & \\
  0 & & K_{71} & & & K_{72} & & & K_{73} & & \ldots & & K_{77} & \\
  0 & & & & & & & & & & & & &
\end{array}
\right)
\end{equation}

For the desired matrix $A$, the conditions~\eqref{def:main1} and~\eqref{def:main2} will take the following form:
\begin{equation}
    \label{eq:main1dsrg}
    A^2 + A = 2I_{22} + 4J_{22},
\end{equation}
\begin{equation}
    \label{eq:main2dsrg}
    A J_{22} = J_{22} A = 9 J_{22}.
\end{equation}

Let us write the matrix $A$ in the form
\begin{equation}
A =
\begin{pmatrix}
    0 & B \\
    D & C
\end{pmatrix},
\end{equation}

\noindent where $C$ is the submatrix of matrix $A$ obtained by deleting the first row and the first column, $B = (11 \ldots 10 \ldots 0)$ is the first row of matrix $A$ without the element $A_{11}$, and $D = (11 \ldots 10001110 \ldots 0)^T$ is the first column of matrix $A$ without the element $A_{11}$.

With this representation of the matrix $A$, the equality~\eqref{eq:main1dsrg} is equivalent to the system of the following four equations:
\begin{equation}
\label{eq:sysBCD}
\begin{cases}
BD = 6,\\
BC + B = 4J_{1,21},\\
CD + D = 4J_{21,1},\\
DB + C^2 + C = 2I_{21} + 4J_{21}.
\end{cases}
\end{equation}

The first equality in the system~\eqref{eq:sysBCD} is satisfied automatically, this follows from the form of the matrices $B$ and $D$.

From the second equality in~\eqref{eq:sysBCD} it follows that
\begin{equation}
\label{eq:sumC1}
\sum_{i=1}^9 C_{ij} =
\begin{cases}
3, \; \text{if} \; j \leqslant 9,\\
4, \; \text{if} \; j \geqslant 10.
\end{cases}
\end{equation}

Since the desired matrix satisfies the relation $AJ_{22} = 9J_{22}$, then from~\eqref{eq:sumC1} it follows that
\begin{equation}
\label{eq:sumC11}
\sum_{i=10}^{21} C_{ij} = 5 \; \text{for} \; 1 \leqslant j \leqslant 21.
\end{equation}

From the third equality in~\eqref{eq:sysBCD} it follows that
\begin{equation}
\label{eq:sumC2}
\sum_{j \in [12]\setminus \{7, 8, 9\}} C_{ij} =
\begin{cases}
3, \; \text{if} \; i \in [12]\setminus \{7, 8, 9\},\\
4, \; \text{if} \; i \notin [12]\setminus \{7, 8, 9\}.
\end{cases}
\end{equation} 

From the equalities~\eqref{eq:main2dsrg} and~\eqref{eq:sumC2} it follows that
\begin{equation}
\label{eq:sumC21}
\sum_{j \notin [12]\setminus \{7, 8, 9\}} C_{ij} = 5 \; \text{for} \; 1 \leqslant i \leqslant 21.
\end{equation}

Let us consider the matrix $H_7$, which has the following form:
\begin{equation}
H_7 =
\begin{pmatrix}
3 & 3 & 3 & 4 & 4 & 4 & 4 \\
3 & 3 & 3 & 4 & 4 & 4 & 4 \\
4 & 4 & 4 & 4 & 4 & 4 & 4 \\
3 & 3 & 3 & 4 & 4 & 4 & 4 \\
4 & 4 & 4 & 4 & 4 & 4 & 4 \\
4 & 4 & 4 & 4 & 4 & 4 & 4 \\
4 & 4 & 4 & 4 & 4 & 4 & 4
\end{pmatrix}.
\end{equation}

We introduce the notation: $H_{21} = 4J_{21} - DB$. By direct computation, one can verify that the matrix $H_{21}$ is obtained from the matrix $H_7$ by replacing every number 3 with blocks $3J_3$ and every number 4 with blocks $4J_3$.

The last equation of the system~\eqref{eq:sysBCD} can be rewritten as
\begin{equation}
\label{eq:small_eq}
C^2 + C = 2I_{21} + H_{21}.
\end{equation} 

We apply the method described in the introduction. Since the matrix $C$ consists of $7 \times 7$ circulant blocks of order three, it can be compactified into a matrix $C(x)$ with elements from $\mathbb{Z}[x]/(x^3 - 1)$.

The equation~\eqref{eq:small_eq} is transformed into the following:
\begin{equation}
\label{eq:small_eq_pol}
C^2(x) + C(x) \equiv 2I_7 + (1+x+x^2)H_7 \pmod{x^3-1}.
\end{equation}

Since $x = 1$ is a root of the polynomial $x^3 - 1$, then when substituting $x = 1$ into~\eqref{eq:small_eq_pol} we obtain the following true equality:
\begin{equation}
\label{eq:small_eq_C1}
C^2(1) + C(1) = 2I_7 + 3H_7.
\end{equation}

Since all coefficients of the polynomial $C(x)$ are zeros or ones, the elements of the matrix $C(1)$ are integers from 0 to 3. The diagonal elements of the matrix $C(1)$ are at most two because the adjacency matrix of the digraph we seek has zeros on the diagonal.

The authors developed a computer program that searches for all seventh-order matrices with elements from 0 to 3 (diagonal elements from 0 to 2) satisfying condition~\eqref{eq:small_eq_C1} and the conditions~\eqref{eq:sumC1}, \eqref{eq:sumC11}, \eqref{eq:sumC2}, \eqref{eq:sumC21}, rewritten for the compactified matrix. The search was implemented using the constraint programming library Artelys Kalis~\cite{bibKalis1}.

As a result of the search, 10338 matrices suitable for the role of the matrix $C(1)$ were obtained. Note that having the matrix $C(1)$ significantly narrows the search space for the adjacency matrix $A$ of the digraph $\text{dsrg}(22, 9, 6, 3, 4)$: restrictions appear on the sums of elements in the rows of the circulant blocks of matrix $A$. Based on these considerations, the list of matrices $C(1)$ was reduced to 144 matrices for which a matrix $A$ satisfying conditions~\eqref{eq:main1dsrg} and~\eqref{eq:main2dsrg} can be found.

For each of these 144 matrices $C(1)$, all suitable adjacency matrices $A$ were found programmatically. There were 384 such matrices. The resulting digraphs will be analyzed in the next section.

The program was run on a computer with an Intel Core i5-7400 processor (3.00 GHz) and 32~GB of RAM. The distribution of the search stages by time is as follows:
\begin{enumerate}
    \item search for 10338 matrices suitable as the matrix $C(1)$: 125 hours;
    \item search for matrices $A$ for the found matrices $C(1)$: 40 hours.
\end{enumerate}
\begin{remark}
    It is worth noting that using the developed program, not all possible adjacency matrices of the digraph $\text{dsrg}(22, 9, 6, 3, 4)$ whose main part consists of $3 \times 3$ circulants were found, but only those whose first row and first column have the specified form (see~\eqref{eq:block_repr}).
\end{remark}

\section{Analysis of constructed examples for isomorphism}

Let us give as an example one of the adjacency matrices found using the computer program:

\begin{equation}
\label{eq:adj_example}
\left(
\begin{array}{cccccccccccccccccccccc}
  0  &  1  &  1  &  1  &  1  &  1  &  1  &  1  &  1  &  1  &  0  &  0  &  0  &  0  &  0  &  0  &  0  &  0  &  0  &  0  &  0  &  0 \\
  1  &  0  &  1  &  1  &  0  &  0  &  0  &  0  &  0  &  1  &  1  &  0  &  0  &  0  &  0  &  0  &  0  &  1  &  0  &  1  &  1  &  1 \\
  1  &  1  &  0  &  1  &  0  &  0  &  0  &  1  &  0  &  0  &  0  &  1  &  0  &  0  &  0  &  0  &  0  &  0  &  1  &  1  &  1  &  1 \\
  1  &  1  &  1  &  0  &  0  &  0  &  0  &  0  &  1  &  0  &  0  &  0  &  1  &  0  &  0  &  0  &  1  &  0  &  0  &  1  &  1  &  1 \\
  1  &  0  &  0  &  1  &  0  &  0  &  0  &  1  &  1  &  0  &  0  &  1  &  1  &  0  &  0  &  1  &  1  &  0  &  0  &  1  &  0  &  0 \\
  1  &  1  &  0  &  0  &  0  &  0  &  0  &  0  &  1  &  1  &  1  &  0  &  1  &  1  &  0  &  0  &  0  &  1  &  0  &  0  &  1  &  0 \\
  1  &  0  &  1  &  0  &  0  &  0  &  0  &  1  &  0  &  1  &  1  &  1  &  0  &  0  &  1  &  0  &  0  &  0  &  1  &  0  &  0  &  1 \\
  0  &  0  &  0  &  0  &  1  &  1  &  1  &  0  &  0  &  0  &  0  &  1  &  0  &  1  &  1  &  1  &  1  &  1  &  0  &  0  &  0  &  0 \\
  0  &  0  &  0  &  0  &  1  &  1  &  1  &  0  &  0  &  0  &  0  &  0  &  1  &  1  &  1  &  1  &  0  &  1  &  1  &  0  &  0  &  0 \\
  0  &  0  &  0  &  0  &  1  &  1  &  1  &  0  &  0  &  0  &  1  &  0  &  0  &  1  &  1  &  1  &  1  &  0  &  1  &  0  &  0  &  0 \\
  1  &  1  &  0  &  0  &  0  &  0  &  0  &  0  &  0  &  1  &  0  &  1  &  1  &  0  &  0  &  0  &  0  &  1  &  0  &  1  &  1  &  1 \\
  1  &  0  &  1  &  0  &  0  &  0  &  0  &  1  &  0  &  0  &  1  &  0  &  1  &  0  &  0  &  0  &  0  &  0  &  1  &  1  &  1  &  1 \\
  1  &  0  &  0  &  1  &  0  &  0  &  0  &  0  &  1  &  0  &  1  &  1  &  0  &  0  &  0  &  0  &  1  &  0  &  0  &  1  &  1  &  1 \\
  0  &  0  &  1  &  1  &  0  &  1  &  0  &  1  &  0  &  0  &  0  &  1  &  0  &  0  &  0  &  0  &  1  &  0  &  1  &  1  &  0  &  1 \\
  0  &  1  &  0  &  1  &  0  &  0  &  1  &  0  &  1  &  0  &  0  &  0  &  1  &  0  &  0  &  0  &  1  &  1  &  0  &  1  &  1  &  0 \\
  0  &  1  &  1  &  0  &  1  &  0  &  0  &  0  &  0  &  1  &  1  &  0  &  0  &  0  &  0  &  0  &  0  &  1  &  1  &  0  &  1  &  1 \\
  0  &  0  &  0  &  1  &  1  &  1  &  1  &  1  &  0  &  1  &  0  &  0  &  0  &  1  &  1  &  1  &  0  &  0  &  0  &  0  &  0  &  0 \\
  0  &  1  &  0  &  0  &  1  &  1  &  1  &  1  &  1  &  0  &  0  &  0  &  0  &  1  &  1  &  1  &  0  &  0  &  0  &  0  &  0  &  0 \\
  0  &  0  &  1  &  0  &  1  &  1  &  1  &  0  &  1  &  1  &  0  &  0  &  0  &  1  &  1  &  1  &  0  &  0  &  0  &  0  &  0  &  0 \\
  0  &  1  &  0  &  0  &  1  &  0  &  0  &  0  &  0  &  1  &  1  &  1  &  0  &  1  &  1  &  0  &  0  &  1  &  1  &  0  &  0  &  0 \\
  0  &  0  &  1  &  0  &  0  &  1  &  0  &  1  &  0  &  0  &  0  &  1  &  1  &  0  &  1  &  1  &  1  &  0  &  1  &  0  &  0  &  0 \\
  0  &  0  &  0  &  1  &  0  &  0  &  1  &  0  &  1  &  0  &  1  &  0  &  1  &  1  &  0  &  1  &  1  &  1  &  0  &  0  &  0  &  0
\end{array}
\right).
\end{equation} 

The developed program makes it possible to construct an entire family of digraphs of the required type. Among the digraphs obtained, there will be a considerable number of isomorphic pairs, and determining the number of constructed digraphs up to isomorphism is not entirely straightforward. Exhaustively enumerating all possible bijections of the vertex sets in search of isomorphisms is both impractical and uninformative.

Certain structural features of the constructed digraphs allow this question to be clarified more efficiently.

\begin{definition}
\label{def:main3}
\begin{enumerate}
\item By construction, the automorphism group of each of the digraphs found contains (with appropriate numbering of the vertices) the permutation: $$(1)(2,3,4)(5,6,7)(8,9,10)(11,12,13)(14,15,16)(17,18,19)(20,21,22).$$ In fact, in all constructed examples, the automorphism group of the digraph coincides with the group $\mathbb{Z}_3$ generated by this permutation.
\item Accordingly, the vertex set of the digraph is naturally partitioned into eight subsets: one singleton (the corresponding vertex will be called \textit{special}) and seven subsets of three elements each (these will be called \textit{floors}).
\item The set of edges leading from one floor to another is invariant under the aforementioned permutation group.
\item One of the types of edge sets allowed by the previous item is the \textit{complete transition}, where from each vertex of one floor, there is exactly one edge to each vertex of another floor.
\item We also call it a complete transition if, instead of one of the floors, there is a singleton set consisting of the special vertex (in this case, there are three edges instead of nine).
\item Additionally, note that \textit{within} a floor, there may be a set of three or six edges of the digraph (a cycle or two opposing cycles), but the set of edges within a floor may also be empty.
\item The \textit{skeleton} of the digraph is defined as a digraph with eight vertices corresponding to the vertex subsets, seven of which may be colored \textit{with four colors} corresponding to the presence or absence of internal edges within a floor. The edges of this digraph are (if present) the complete transitions described above.
\item A digraph with the same vertices and (informally speaking) all the remaining edges of the considered digraph (there are several variants of \textit{non-complete} transitions from one floor to another, which are treated as different edge colors) will be called the \textit{rigging}.
\end{enumerate}    
\end{definition}

\begin{figure}[!ht]
    \includegraphics[width=0.49\textwidth]{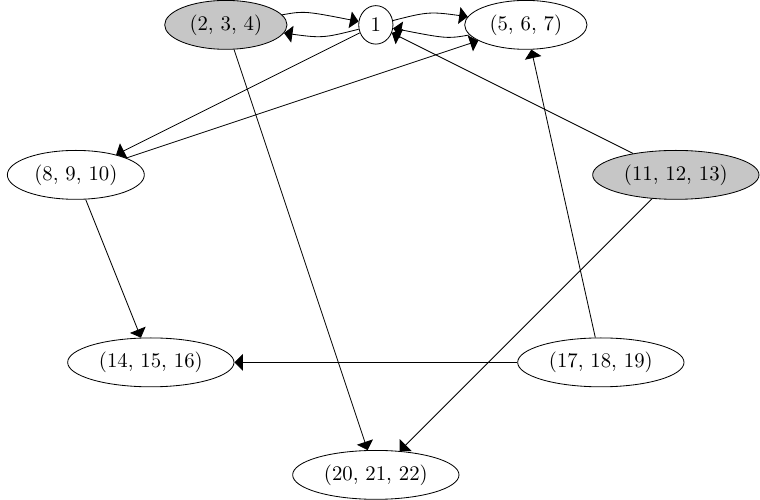}\hfill
    \includegraphics[width=0.49\textwidth]{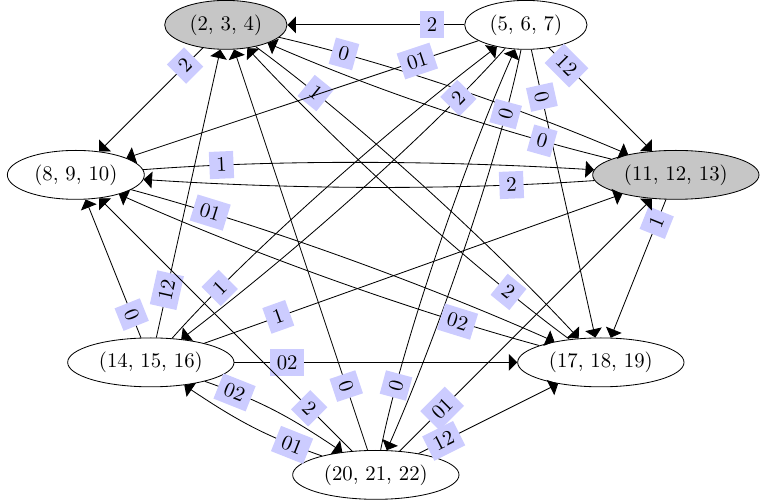}
    \caption{Skeleton and rigging of the first digraph}
    \label{fig:graph1}
\end{figure}

\begin{figure}[!ht]
    \includegraphics[width=0.49\textwidth]{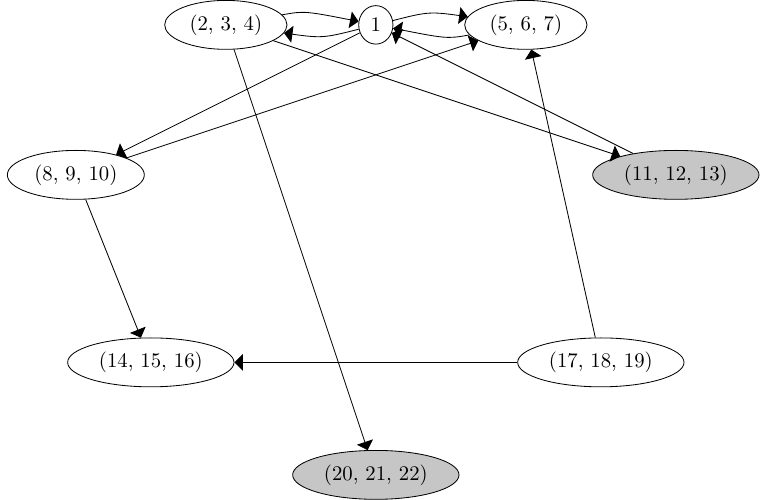}\hfill
    \includegraphics[width=0.49\textwidth]{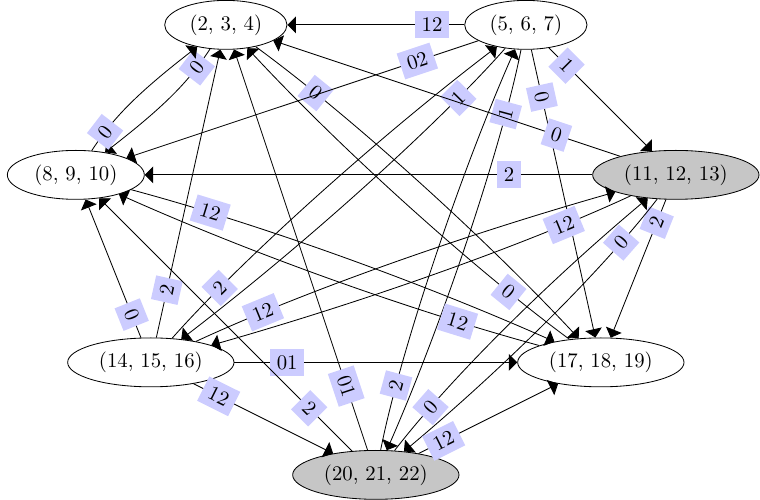}
    \caption{Skeleton and rigging of the second digraph}
    \label{fig:graph2}
\end{figure}

\begin{figure}[!ht]
    \includegraphics[width=0.49\textwidth]{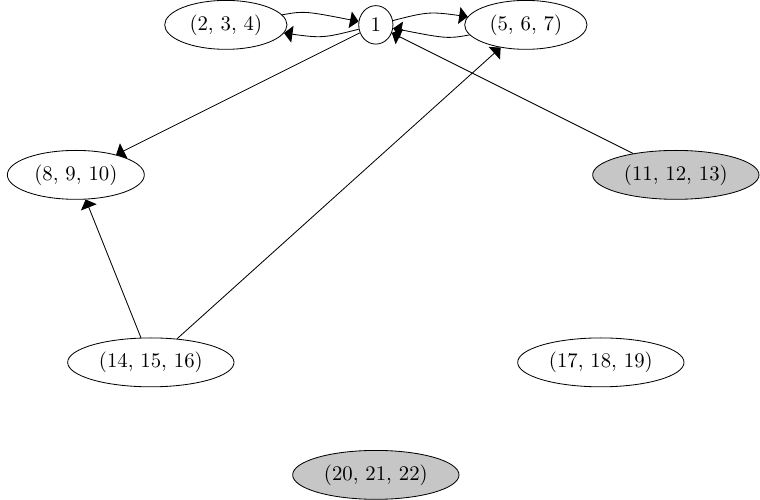}\hfill
    \includegraphics[width=0.49\textwidth]{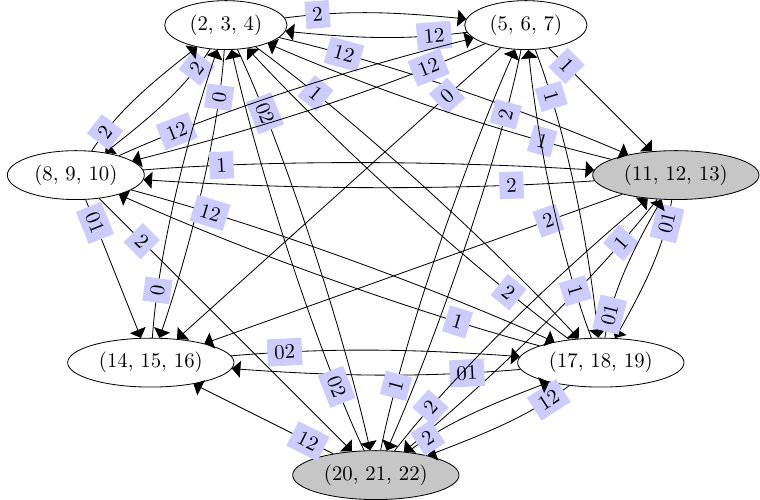}
    \caption{Skeleton and rigging of the third digraph}
    \label{fig:graph3}
\end{figure}

\begin{figure}[!ht]
    \includegraphics[width=0.49\textwidth]{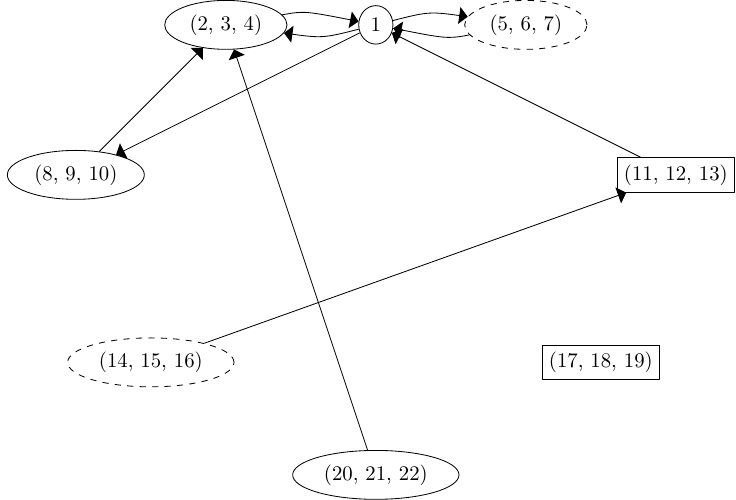}\hfill
    \includegraphics[width=0.49\textwidth]{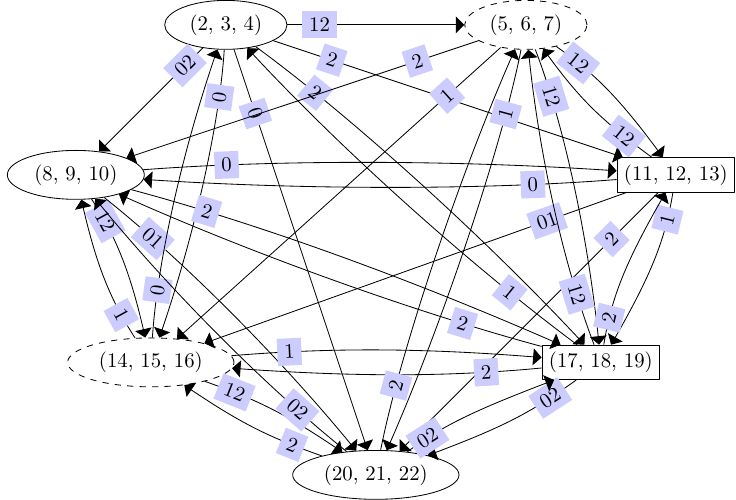}
    \caption{Skeleton and rigging of the fourth digraph}
    \label{fig:graph4}
\end{figure}

\begin{figure}[!ht]
    \includegraphics[width=0.49\textwidth]{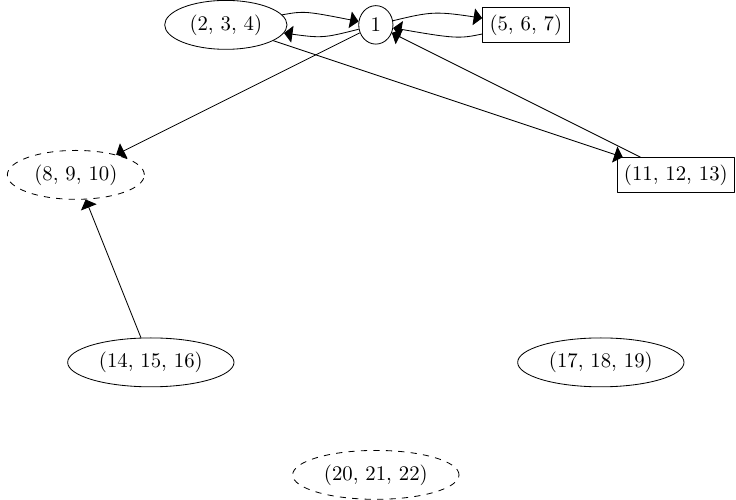}\hfill
    \includegraphics[width=0.49\textwidth]{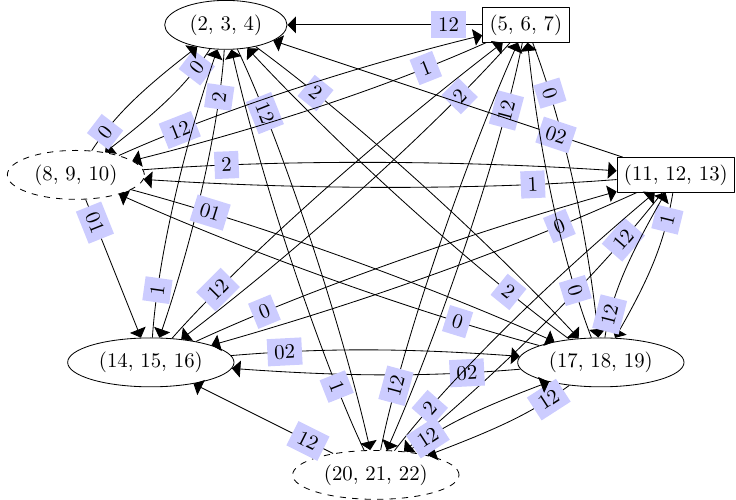}
    \caption{Skeleton and rigging of the fifth digraph}
    \label{fig:graph5}
\end{figure}

\begin{figure}[!ht]
    \includegraphics[width=0.49\textwidth]{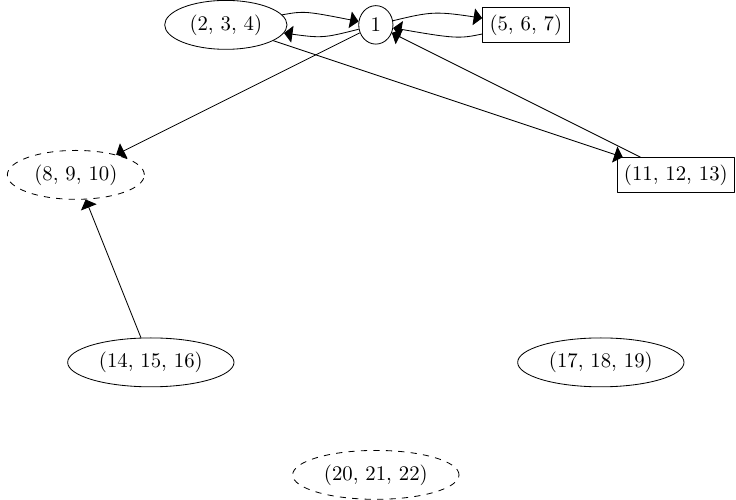}\hfill
    \includegraphics[width=0.49\textwidth]{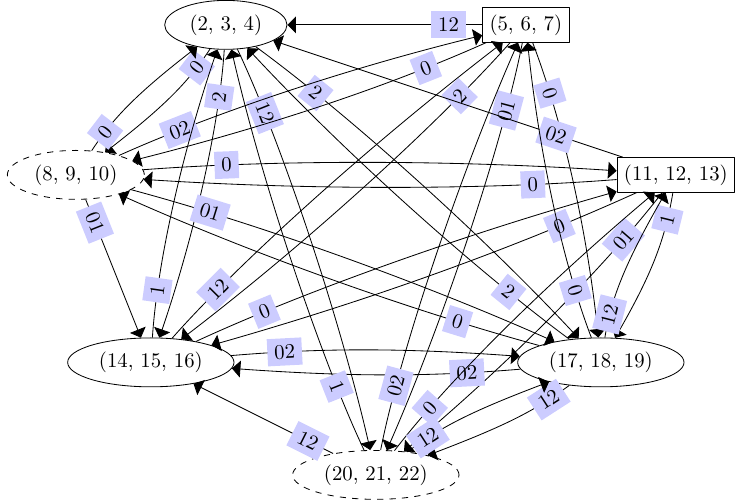}
    \caption{Skeleton and rigging of the sixth digraph}
    \label{fig:graph6}
\end{figure}

\begin{figure}[!ht]
    \includegraphics[width=0.49\textwidth]{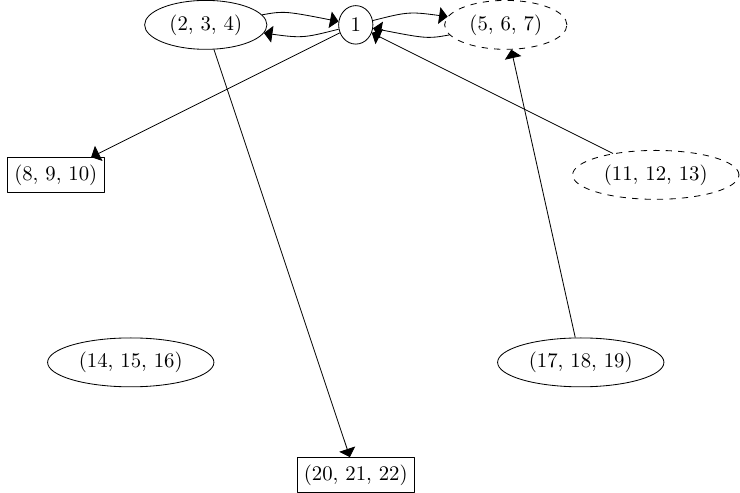}\hfill
    \includegraphics[width=0.49\textwidth]{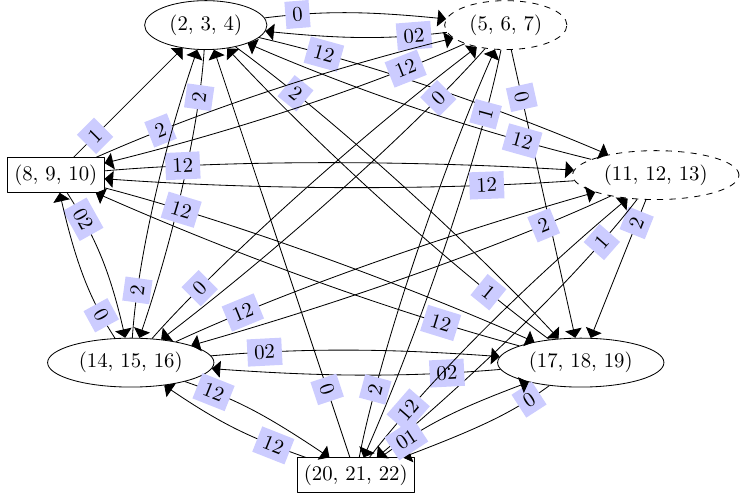}
    \caption{Skeleton and rigging of the seventh digraph}
    \label{fig:graph7}
\end{figure}

\begin{figure}[!ht]
    \includegraphics[width=0.49\textwidth]{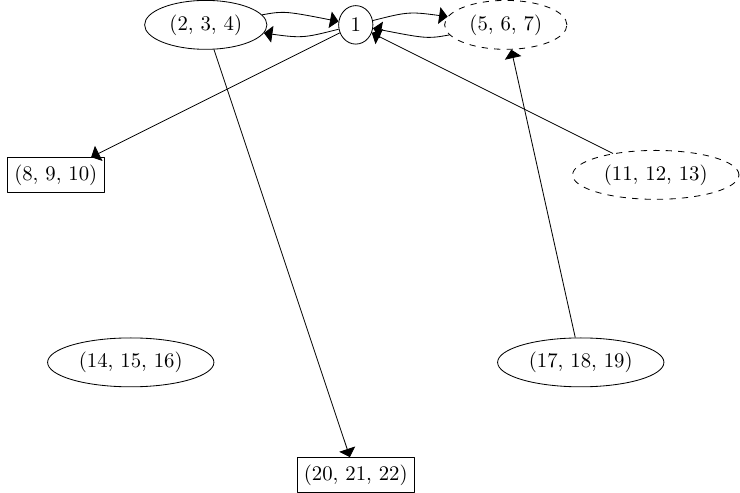}\hfill
    \includegraphics[width=0.49\textwidth]{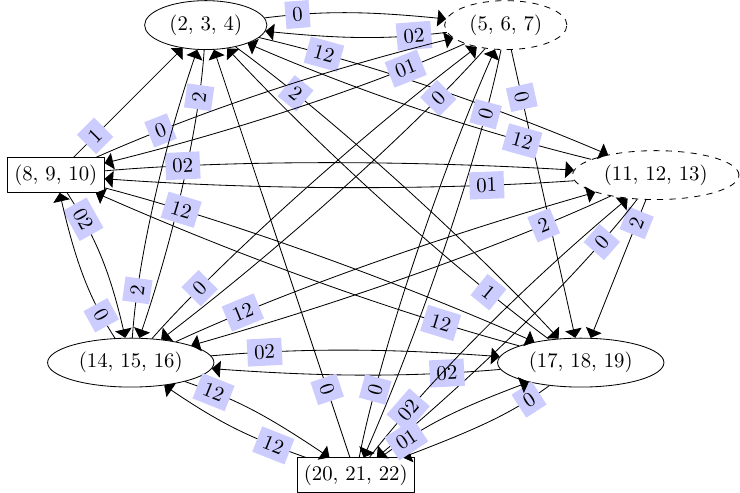}
    \caption{Skeleton and rigging of the eighth digraph}
    \label{fig:graph8}
\end{figure}

The digraphs found using the computer program are divided into 16 groups, each group consisting of 24 isomorphic digraphs. Figures~\ref{fig:graph1}--\ref{fig:graph8} show the skeletons and rigging corresponding to eight groups (in the rigging diagrams, the special vertex is not shown since it is isolated). The remaining eight digraphs are obtained by reversing the direction of all edges in the given digraphs.

The following notations are used in the figures shown.

\begin{enumerate}
\item If within a floor no vertex is connected to another by an edge, then the corresponding floor is depicted as an unfilled ellipse with a solid border (for example, $(5, 6, 7)$ in Fig.~\ref{fig:graph1}).
\item If within a floor every two vertices are connected by a pair of edges, then the corresponding floor is depicted as a filled ellipse (for example, $(2, 3, 4)$ in Fig.~\ref{fig:graph1}).
\item If a floor contains one cycle consisting of edges $(3k+2) \to (3k+3)$, \mbox{$(3k+3) \to (3k+4)$}, \mbox{$(3k+4) \to (3k+2)$} ($k$ takes one of the values: $0, 2, \ldots, 6$), then the corresponding floor is depicted as an ellipse with a dashed border (for example, $(5, 6, 7)$ in Fig.~\ref{fig:graph4}).
\item If a floor contains one cycle consisting of edges $(3k+4) \to (3k+3)$, \mbox{$(3k+3) \to (3k+2)$}, \mbox{$(3k+2) \to (3k+4)$} ($k$ takes one of the values: $0, 2, \ldots, 6$), then the corresponding floor is depicted as a rectangle (for example, $(11, 12, 13)$ in Fig.~\ref{fig:graph4}).
\item If the label on an edge $(a_0, a_1, a_2) \to (b_0, b_1, b_2)$ in the rigging contains a digit $s$ (\mbox{$s \in {0, 1, 2}$}), then in the original digraph there is an edge from each vertex $a_{i}$ to vertex $b_{(i+s) \pmod{3}}$ for $i \in {0, 1, 2}$. For example, in the rigging from Fig.~\ref{fig:graph1}, there is an edge from floor $(14, 15, 16)$ to floor $(2, 3, 4)$ labeled $12$, which means the presence in the original digraph of edges $14\to 3$, $15\to 4$, $16\to 2$ and $14\to 4$, $15\to 2$, $16\to 3$.
\end{enumerate}

For a digraph whose adjacency matrix is given in the formula~\eqref{eq:adj_example}, the skeleton and rigging are shown in Fig.~\ref{fig:graph1}.

Note that the digraph's skeleton does not uniquely define its rigging: the fifth and sixth digraphs (Fig.~\ref{fig:graph5} and ~\ref{fig:graph6}) have the same skeletons, but different riggings. The seventh and eighth digraphs also have the same skeletons (Fig.~\ref{fig:graph7} and ~\ref{fig:graph8}).

\begin{remark}
    It should be noted that the question of investigating the found digraphs for isomorphism was also solved programmatically using the is\_isomorphic function from the NetworkX library~\cite{bibNetworkx1}. However, the authors believe that theoretical analysis of graph structure is more informative than simply stating their count using a program. The main goal of the theoretical analysis was to construct a system of invariants, analogous to which may be encountered in other problems.
\end{remark}

\begin{remark}
    In \cite{bibDuval1}, the following fact is proved: if a digraph $G$ is a strongly regular digraph with parameter set $(v, k, t, \lambda, \mu)$, then the digraph $G^{\prime}$, which is the complement of $G$, is also a strongly regular digraph with parameter set $(v, v-k-1, v-2k+t-1, v-2k+\mu-2, v-2k+\lambda)$ (the adjacency matrix of graph $G^{\prime}$ equals $A'=J_v-I_v-A$, where $A$ is the adjacency matrix of graph $G$). It follows that this work has also automatically constructed a family of strongly regular digraphs with parameter set $(22, 12, 9, 6, 7)$.
\end{remark}

\section{Conclusion}

In this work, we have constructed a family of strongly regular digraphs with parameter set $(22, 9, 6, 3, 4)$, the existence of which was previously considered an open question. The construction of the matrix is based on an idea of O.~Gritsenko presented in \cite{bibGritsenko1}. In addition, we have proposed a natural structure that allows the question of their isomorphism to be resolved in a simple, intuitively clear manner (without having to rely on exhaustive search programs). The structural decomposition proved to be nontrivial in the sense that the first structural element (skeleton) does not always determine the second (rigging) uniquely, and useful in the sense that it (theoretically) allows manual verification that a digraph satisfies the required constraints.

Obviously, the described construction (partitioning the part of the matrix without the first row and first column into circulant blocks) is applicable only to parameter sets $(v, k, t, \lambda, \mu)$ for which the numbers $v-1$, $k$, $t$ are divisible by the block size. This situation frequently occurs among problems that are currently unsolved, so the method used by the authors is potentially applicable to them as well; however, the authors are unaware of other advances of this type.

\bibliographystyle{unsrtnat}
%\bibliography{references}  %%% Uncomment this line and comment out the ``thebibliography'' section below to use the external .bib file (using bibtex) .


\end{document}